# ON CARMICHAËL'S CONJECTURE

Florentin Smarandache
University of New Mexico
200 College Road
Gallup, NM 87301, USA
E-mail: smarand@unm.edu

**Introduction**.
Carmichaël's conjecture is the following: "the equation $\varphi(x) = n$ cannot have a unique solution, $(\forall) n \in \mathbb{N}$, where $\varphi$ is the Euler's function". R. K. Guy presented in [1] some results on this conjecture; Carmichaël himself proved that, if $n_0$ does not verify his conjecture, then $n_0 > 10^{37}$; V. L. Klee [2] improved to $n_0 > 10^{400}$, and P. Masai & A. Valette increased to $n_0 > 10^{10000}$. C. Pomerance [4] wrote on this subject too.

In this article we prove that the equation $\varphi(x) = n$ admits a finite number of solutions, we find the general form of these solutions, also we prove that, if $x_0$ is the unique solution of this equation (for a $n \in \mathbb{N}$), then $x_0$ is a multiple of $2^2 \cdot 3^2 \cdot 7^2 \cdot 43^2$ (and $x_0 > 10^{10000}$ from [3]).

In the last paragraph we extend the result to: $x_0$ is a multiple of a product of a very large number of primes.

§1. Let $x_0$ be a solution of the equation $\varphi(x) = n$. We consider $n$ fixed. We'll try to construct another solution $y_0 \neq x_0$.

*The first method:*

We decompose $x_0 = a \cdot b$ with $a$, $b$ integers such that $(a, b) = 1$.

we look for an $a' \neq a$ such that $\varphi(a') = \varphi(a)$ and $(a', b) = 1$; it results that $y_0 = a' \cdot b$.

*The second method:*

Let's consider $x_0 = q_1^{\beta_1} ... q_r^{\beta_r}$, where all $\beta_i \in \mathbb{N}^*$, and $q_1, ..., q_r$ are distinct primes two by two; we look for an integer $q$ such that $(q, x_0) = 1$ and $\varphi(q)$ divides $x_0 / (q_1, ..., q_r)$; then $y_0 = x_0 q / \varphi(q)$.

We immediately see that we can consider $q$ as prime.

The author conjectures that for any integer $x_0 \geq 2$ it is possible to find, by means of one of these methods, a $y_0 \neq x_0$ such that $\varphi(y_0) = \varphi(x_0)$.

**Lemma 1.** The equation $\varphi(x) = n$ admits a finite number of solutions, $(\forall) n \in \mathbb{N}$.
*Proof.* The cases $n = 0, 1$ are trivial.

Let's consider $n$ to be fixed, $n \geq 2$. Let $p_1 < p_2 < ... < p_s \leq n+1$ be the sequence of prime numbers. If $x_0$ is a solution of our equation (1) then $x_0$ has the form $x_0 = p_1^{\alpha_1}...p_s^{\alpha_s}$, with all $\alpha_i \in \mathbb{N}$. Each $\alpha_i$ is limited, because:
$$(\forall) i \in \{1,2,...,s\}, \ (\exists) a_i \in \mathbb{N}: \ p_i^{a_i} \geq n.$$
Whence $0 \leq \alpha_i \leq a_i + 1$, for all $i$. Thus, we find a wide limitation for the number of solutions: $\prod_{i=1}^{s}(a_i + 2)$

**Lemma 2.** Any solution of this equation has the form (1) and (2):
$$x_0 = n \cdot \left(\frac{p_1}{p_1 - 1}\right)^{\varepsilon_1} \cdot ... \cdot \left(\frac{p_s}{p_s - 1}\right)^{\varepsilon_s} \in \mathbb{Z},$$
where, for $1 \leq i \leq s$, we have $\varepsilon_i = 0$ if $\alpha_i = 0$, or $\varepsilon_i = 1$ if $\alpha_i \neq 0$.

Of course, $n = \varphi(x_0) = x_0 \left(\frac{p_1}{p_1 - 1}\right)^{\varepsilon_1} \cdot ... \cdot \left(\frac{p_s}{p_s - 1}\right)^{\varepsilon_s}$,

whence it results the second form of $x_0$.

From (2) we find another limitation for the number of the solutions: $2^s - 1$ because each $\varepsilon_i$ has only two values, and at least one is not equal to zero.

§2. We suppose that $x_0$ is the unique solution of this equation.

**Lemma 3.** $x_0$ is a multiple of $2^2 \cdot 3^2 \cdot 7^2 \cdot 43^2$.

*Proof.* We apply our second method.

Because $\varphi(0) = \varphi(3)$ and $\varphi(1) = \varphi(2)$ we take $x_0 \geq 4$.

If $2 \nmid x_0$ then there is $y_0 = 2x_0 \neq x_0$ such that $\varphi(y_0) = \varphi(x_0)$, hence $2 \mid x_0$; if $4 \nmid x_0$, then we can take $y_0 = x_0 / 2$.

If $3 \nmid x_0$ then $y_0 = 3x_0 / 2$, hence $3 \mid x_0$; if $9 \nmid x_0$ then $y_0 = 2x_0 / 3$, hence $9 \mid x_0$; whence $4 \cdot 9 \mid x_0$.

If $7 \nmid x_0$ then $y_0 = 7x_0 / 6$, hence $7 \mid x_0$; if $49 \nmid x_0$ then $y_0 = 6x_0 / 7$ hence $49 \mid x_0$; whence $4 \cdot 9 \cdot 49 \mid x_0$.

If $43 \nmid x_0$ then $y_0 = 43x_0 / 42$, hence $43 \mid x_0$; if $43^2 \nmid x_0$ then $y_0 = 42x_0 / 43$, hence $43^2 \mid x_0$; whence $2^2 \cdot 3^2 \cdot 7^2 \cdot 43^2 \mid x_0$.

Thus $x_0 = 2^{\gamma_1} \cdot 3^{\gamma_2} \cdot 7^{\gamma_3} \cdot 43^{\gamma_4} \cdot t$, with all $\gamma_i \geq 2$ and $(t, \ 2 \cdot 3 \cdot 7 \cdot 43) = 1$ and $x_0 > 10^{10000}$ because $n_0 > 10^{10000}$.

§3. Let's consider $\gamma_1 \geq 3$. If $5 \nmid x_0$ then $5x_0 / 4 = y_0$, hence $5 \mid x_0$; if $25 \nmid x_0$ then $y_0 = 4x_0 / 5$, whence $25 \mid x_0$.

We construct the recurrent set $M$ of prime numbers:
a) the elements $2, 3, 5 \in M$;
b) if the distinct odd elements $e_1,...,e_n \in M$ and $b_m = 1 + 2^m \cdot e_1,...,e_n$ is prime, with $m = 1$ or $m = 2$, then $b_m \in M$;

c) any element belonging to $M$ is obtained by the utilization (a finite number of times) of the rules a) or b) only.

The author conjectures that $M$ is infinite, which solves this case, because it results that there is an infinite number of primes which divide $x_0$. This is absurd.

For example 2, 3, 5, 7, 11, 13, 23, 29, 31, 43, 47, 53, 61, … belong to $M$.

\*

The method from §3 could be continued as a tree (for $\gamma_2 \geq 3$ afterwards $\gamma_3 \geq 3$, etc.) but its ramifications are very complicated…

§4. **A Property for a Counter-Example to Carmichael Conjecture**.
Carmichaël has conjectured that:
$(\forall)\, n \in \mathbb{N}$, $(\exists)\, m \in \mathbb{N}$, with $m \neq n$, for which $\varphi(n) = \varphi(m)$, where $\varphi$ is Euler's totient function.

There are many papers on this subject, but the author cites the papers which have influenced him, especially Klee's papers.

Let n be a counterexample to Carmichaël's conjecture.

Grosswald has proved that $n_0$ is a multiple of 32, Donnelly has pushed the result to a multiple of $2^{14}$, and Klee to a multiple of $2^{42} \cdot 3^{47}$, Smarandache has shown that n is a multiple of $2^2 \cdot 3^2 \cdot 7^2 \cdot 43^2$. Masai & Valette have bounded
$n > 10^{10000}$.

In this paragraph we will extend these results to: $n$ is a multiple of a product of a very large number of primes.

We construct a recurrent set $M$ such that:
a) the elements $2, 3 \in M$;
b) if the distinct elements $2, 3, q_1, ..., q_r \in M$ and $p = 1 + 2^a \cdot 3^b \cdot q_1 \cdots q_r$ is a prime, where $a \in \{0, 1, 2, ..., 41\}$ and $b \in \{0, 1, 2, ..., 46\}$, then $p \in M$; $r \geq 0$;
c) any element belonging to $M$ is obtained only by the utilization (a finite number of times) of the rules a) or b).

Of course, all elements from $M$ are primes.

Let $n$ be a multiple of $2^{42} \cdot 3^{47}$;
if $5 \nmid n$ then there exists $m = 5n/4 \neq n$ such that $\varphi(n) = \varphi(m)$; hence
$5 \mid n$; whence $5 \in M$;
if $5^2 \nmid n$ then there exists $m = 4n/5 \neq n$ with our property; hence $5^2 \mid n$;
analogously, if $7 \nmid n$ we can take $m = 7n/6 \neq n$, hence $7 \mid n$; if $7^2 \nmid n$ we can take $m = 6n/7 \neq n$; whence $7 \in M$ and $7^2 \mid n$; etc.

The method continues until it isn't possible to add any other prime to $M$, by its construction.

For example, from the 168 primes smaller than 1000, only 17 of them do not belong to $M$ (namely: 101, 151, 197, 251, 401, 491, 503, 601, 607, 677, 701, 727, 751, 809, 883, 907, 983); all other 151 primes belong to $M$.

Note $M = \{2, 3, p_1, p_2, ..., p_s, ...\}$, then $n$ is a multiple of $2^{42} \cdot 3^{47} \cdot p_1^2 \cdot p_2^2 \cdots p_s^2 \cdots$
From our example, it results that $M$ contains at least 151 elements, hence $s \geq 149$.

If $M$ is infinite then there is no counterexample $n$, whence Carmichaël's conjecture is solved.

(The author conjectures $M$ is infinite.)

Using a computer it is possible to find a very large number of primes, which divide $n$, using the construction method of $M$, and trying to find a new prime $p$ if $p - 1$ is a product of primes only from $M$.